\documentclass[12pt,a4paper]{article}
\usepackage[cp1251]{inputenc}
\usepackage[T2A]{fontenc}
\usepackage[english]{babel}

\usepackage{amsmath,amsfonts,amssymb,mathrsfs,amsopn,indentfirst,amscd}

\usepackage{graphicx}
\usepackage{amsthm}
\usepackage{hyperref}

\usepackage{amssymb}
\usepackage{amsmath}

\newcommand{\vo}{{\rm vol}}

\newcommand{\conv}{{\rm conv}}

\newtheorem*{corollary*}{Corollary}

\begin{document}

\title{On Properties of a Regular Simplex  \\
Inscribed into a Ball}
\author{Mikhail Nevskii\footnote{ Department of Mathematics,  P.G.~Demidov Yaroslavl State University,\newline Sovetskaya str., 14, Yaroslavl, 150003, Russia, \newline
               mnevsk55@yandex.ru,  orcid.org/0000-0002-6392-7618 } 
               }
      
\date{May 17, 2021}
\maketitle

\begin{abstract}

\medskip

  Let  $B$ be a Euclidean ball in ${\mathbb R}^n$ and let 
$C(B)$ be a space of~continuous functions 
$f:B\to{\mathbb R}$ with the uniform norm 
$\|f\|_{C(B)}:=\max_{x\in B}|f(x)|.$ By~$\Pi_1\left({\mathbb R}^n\right)$ we mean a set of polynomials of degree
$\leq 1$, i.\,e., a set of linear functions upon
${\mathbb R}^n$. The interpolation projector  
$P:C(B)\to \Pi_1({\mathbb R}^n)$ with the nodes $x^{(j)}\in B$
is defined by~the~equalities
$Pf\left(x^{(j)}\right)=
f\left(x^{(j)}\right)$,  $j=1,$ $\ldots,$ $ n+1$.
The norm of $P$ as an operator from
$C(B)$ to $C(B)$ can be calculated  by the formula 
$\|P\|_B=\max_{x\in B}\sum |\lambda_j(x)|.$ Here $\lambda_j$ are the basic Lagrange
polynomials corresponding to the $n$-dimensional nondegenerate simplex $S$ 
with the vertices  $x^{(j)}$. Let $P^\prime$ be a projector having the nodes in the vertices \linebreak of a regular simplex inscribed into the ball. 
We describe the points $y\in B$ with the property $\|P^\prime\|_B=\sum |\lambda_j(y)|$. 
Also we formulate a geometric conjecture which implies that $\|P^\prime\|_B$ is equal to the
minimal norm of an interpolation projector with nodes in $B$.  We prove that this conjecture holds true at least for~$n=1,2,3,4$. 

\medskip 


\end{abstract}


\section{Introduction}\label{nev_s1}

Let $\Omega$ be a convex body in
${\mathbb R}^n$. Denote by   
$C(\Omega)$ a space of continuous functions 
$f:\Omega\to{\mathbb R}$ with  the uniform norm
$$\|f\|_{C(\Omega)}:=\max\limits_{x\in Q_n}|f(x)|.$$
By $\Pi_1\left({\mathbb R}^n\right)$ we mean a set of polynomials in
$n$ variables of degree
$\leq 1$, i.\,e.,  of  linear functions on ${\mathbb R}^n$.
For $x^{(0)}\in {\mathbb R}^n,  R>0$, by $B(x^{(0)};R)$ we denote the $n$-dimensional Euclidean ball given
by the inequality 
$\|x-x^{(0)}\|\leq R$. Here 
$$\|x\|:=\sqrt{(x,x)}=\left(\sum\limits_{i=1}^n x_i^2\right)^{1/2}.$$
By definition,  $B_n:=B(0;1)$. Further $e_1,\ldots,e_n$ is the canonical basis in  ${\mathbb R}^n$.

Let $S$ be a nondegenerate simplex in  ${\mathbb R}^n$ with vertices 
$x^{(j)}=\left(x_1^{(j)},\ldots,x_n^{(j)}\right),$ 
$1\leq j\leq n+1.$ Consider the following {\it vertex matrix} of this simplex:
$${\bf S} :=
\left( \begin{array}{cccc}
x_1^{(1)}&\ldots&x_n^{(1)}&1\\
x_1^{(2)}&\ldots&x_n^{(2)}&1\\
\vdots&\vdots&\vdots&\vdots\\
x_1^{(n+1)}&\ldots&x_n^{(n+1)}&1\\
\end{array}
\right).$$
Let ${\bf S}^{-1}$ $=(l_{ij})$. Linear  polynomials 
$\lambda_j(x):=
l_{1j}x_1+\ldots+
l_{nj}x_n+l_{n+1,j}$
have a~property  
$\lambda_j\left(x^{(k)}\right)$ $=$ 
$\delta_j^k$.
We call $\lambda_j$ {\it the basic Lagrange polynomials corresponding to $S$.}
For~an~arbitrary $x\in{\mathbb R}^n$, 
$$x=\sum_{j=1}^{n+1} \lambda_j(x)x^{(j)}, \quad \sum_{j=1}^{n+1} \lambda_j(x)=1.$$
These equalities mean that $\lambda_j(x)$ are {\it  the barycentric coordinates of $x$}. 
For details, see~\cite[\S 1.1]{nevskii_monograph}.

An interpolation projector
 $P:C(\Omega)\to \Pi_1({\mathbb R}^n)$ corresponds to a simplex $S\subset\Omega$  if the nodes of $P$
 coincide with the vertices of $S$. This projector is defined by the~equalities
$Pf\left(x^{(j)}\right)=
f\left(x^{(j)}\right).$
The following analogue of the Lagrange interpolation formula holds:
\begin{equation}\label{interp_Lagrange_formula}
Pf(x)=\sum\limits_{j=1}^{n+1}
f\left(x^{(j)}\right)\lambda_j(x). 
\end{equation}
Denote by $\|P\|_\Omega$ the norm of $P$ as an operator from $C(\Omega)$
in $C(\Omega)$. From (\ref{interp_Lagrange_formula}), it~follows that
$$
\|P\|_\Omega=
\max_{x\in\Omega}\sum_{j=1}^{n+1}
|\lambda_j(x)|.
$$
Since $\lambda_j(x)$ are  the barycentric coordinates of a point $x$  we have also
\begin{equation}\label{norm_P_Omega_by_barycentric_coord}
\|P\|_\Omega=\max \left \{ \sum_{j=1}^{n+1} |\beta_j|:  \quad \sum_{j=1}^{n+1} \beta_jx^{(j)} \in \Omega, \quad \sum_{j=1}^{n+1} \beta_j=1\right\}.
\end{equation}

Consider the case $\Omega=B:=B(x^{(0)};R)$.   It is proved in \cite{nev_ukh_mais_2019_2} that    
\begin{equation}\label{norm_P_ball_general}
\|P\|_B=
\max\limits_{f_j=\pm 1} \left[
R \left(\sum_{i=1}^n\left(\sum_{j=1}^{n+1} f_jl_{ij}\right)^2\right)^{1/2}
+\left|\sum_{j=1}^{n+1}f_j\lambda_j(x^{(0)})\right|
\right].
\end{equation}
If $S$ is a regular simplex inscribed into the ball, then $\|P\|_B$ depends neither on the 
center $x^{(0)}$ nor on the radius $R$  of the ball nor on the choice of  such a simplex.
In~this case
(see \cite[Theorem 2]{nev_ukh_mais_2019_2})
\begin{equation}\label{norm_P_reg_formula}
\|P\|_B=\max\{\psi(a),\psi(a+1)\}, 
\end{equation}
where $a:=\left\lfloor\frac{n+1}{2}-\frac{\sqrt{n+1}}{2}\right\rfloor,$
\begin{equation}\label{psi_function_modulus}
\psi(t):=\frac{2\sqrt{n}}{n+1}\Bigl(t(n+1-t)\Bigr)^{1/2}+
\left|1-\frac{2t}{n+1}\right|, \quad  0\leq t\leq n+1.
\end{equation} 

Various geometric estimates concerning polynomial interpolation are given in 
\cite{nevskii_monograph}. In particular, this book contains the results corresponding to the linear interpolation on the
unit cube $Q_n:=[0,1]^n$. Later some  estimates for concrete $n$ were  improved 
(e.\,g., see  \cite{nev_ukh_mais_2018_25_3},
\cite{nevskii_ukhalov_2018}).
Interpolation by  linear functions on a Euclidean ball in ${\mathbb R}^n$ and related questions were considered in
\cite{nevskii_mais_2018_6},
\cite{nevskii_mais_2019_3},
\cite{nev_ukh_mais_2019_2}. 

The present paper supplements the results  obtained in \cite{nev_ukh_mais_2019_2}
for a regular simplex inscribed into a ball.
In Section 2, we find the maximal points of the function 
$\lambda(x):=\sum |\lambda_j(x)|$. In any such a point  $y\in B$
$$\lambda(y)=\sum_{j=1}^{n+1} |\lambda_j(y)|=\|P\|_B=\max\{\psi(a),\psi(a+1)\}.$$
The points belong to the boundary sphere of  $B$. The number of these points \linebreak is equal to $N=$ $n+1\choose k$, where
$k$ coincides with that  number $a$ or $a+1$ on which $\psi(t)$ takes a bigger value. 
 Evidently, $N$ is the number of $(k-1)$-dimensional faces of an $n$-dimensional  nondegenerate simplex.
In Section 3 we discuss the questions related to the projector's norm invariance under an affine transform.
In Section 4 we formulate  a geometric conjecture the validity of which implies that the projector corresponding to a regular inscribed simplex has the minimal norm. This  conjecture holds true at least for $n=1,2,3,4$. 

\section{The Maximum Points of $\lambda(x)$  for a Regular \\ Inscribed Simplex}\label{nev_s2}
By definition, put
$$k=k(n):= \left\{
\begin{array}{cc}
a+1,  &  \mbox{���� } \psi(a+1)\geq\psi(a),\\
a,& \mbox{���� } \psi(a+1)< \psi(a).\\
\end{array}\right. 
$$
The detailed analysis and first values $n$ and $k$ are given in   
\cite{nev_ukh_mais_2019_2}. 
For these  data and also the numbers
 $N={n+1\choose k}$, see Table 1. 

For $n=1, 2, 3$, we have  $k=1$.  If $n>3$, then 
$\sqrt{n+1}>2$, hence,
$$a+1=
\left\lfloor\frac{n+1}{2}-\frac{\sqrt{n+1}}{2}\right\rfloor+1\leq
\frac{n+1}{2}-\frac{\sqrt{n+1}}{2}+1<\frac{n+1}{2}.$$
Thus, for $n>3$ holds $k<\frac{n+1}{2}$. Since $k$ is an integer, for all 
$n\geq 2$ we have
$k\leq\frac{n}{2}$.

\begin{table}[h!]
\begin{center}
Table 1. The numbers $n$, $k$, and $N=$ $n+1\choose k$
\bigskip

$
\def\arraystretch{1.7}
\begin{array}{|c|c|c|}
\hline
n & k& N={n+1\choose k} \\
\hline
 1 &1 & 2\\
 \hline
 2 &1 & 3\\
 \hline
 3 &1 & 4\\
 \hline
 4 &1 & 5\\
 \hline
 5 &2 & 15\\
 \hline
 6 &2 & 21\\
 \hline
 7 &3 & 56\\
 \hline
 8 & 3 & 84\\
 \hline
 9 & 3 & 120\\
 \hline
 10 & 4 & 330\\
 \hline
 11 & 4 & 495\\
 \hline
 12 & 5 & 1287\\
 \hline
 13 & 5 & 2002\\
 \hline
 14 & 6 &  5005\\
 \hline
 15 & 6 & 8008\\
 \hline
 50 & 22 & 196793068630200\\
 \hline
 100 & 45 & 110826707011209895344085355160\\
 \hline
\end{array}
$
\end{center}
\end{table}

Suppose $S$ is a  regular simplex inscribed into a ball  $B$,  $\lambda_j$ are the basic Lagrange polynomials of this simplex, and  
$P:C(B)\to \Pi_1({\mathbb R}^n)$ is the corresponding interpolation projector.

 \smallskip
 {\bf Theorem 1.} {\it Consider an arbitrary
$(k-1)$-dimensional face $G$ of the simplex
$S$. Let $H$ be the $(n-k)$-dimensional face of $S$
which contains the vertices not belonging to $G$. Denote by $g$ and $h$ 
the centers of gravity of $G$ and $H$. Assume that $y$ \linebreak  is a point there the 
straight line
$(gh)$ inersects the boundary sphere in direction from
$g$ to $h$. Then
\begin{equation}\label{max_point_y}
\lambda(y)=\sum_{j=1}^{n+1}|\lambda_j(y)|=\|P\|_B.
\end{equation}
}


\smallskip
{\it Proof}. It is sufficient to consider some  given ball $B\subset {\mathbb R}^n$, 
some regular simplex $S$ inscribed into $B$, and also $G=\conv\left(x^{(1)},\ldots,
x^{(k)}\right).$

If $n=1$, then
$\psi(t)=\sqrt{t(2-t)}+|1-t|$, $a=0$, $\psi(a)=\psi(a+1)=1,$ $k=a+1=1$.
Let us take $x^{(1)}=0,$ $x^{(2)}=1,$ i.\,e., $S=B=[0,1].$ In this case, $g=0,$
$h=1$,
$\lambda_1(x)=-x+1,$ $\lambda_2(x)=x$, $\lambda(x) \equiv 1$.
Since $\|P\|_B=1$, for $y=h=1$ the equality (\ref{max_point_y}) holds true.
Note that in this trivial case a set of maximum points of the function $\lambda(x)$ 
coincides with the ball $B$.

Now let $n\geq 2$. 
First note that the ceter of gravity$c$ of the simplex $S$ belongs to~the segment $[g,h]$.
Indeed, the equalities
\begin{equation}\label{c_g_h}
c=\frac{1}{n+1}\sum_{j=1}^{n+1} x^{(j)}, \quad
g=\frac{1}{k}\sum_{j=1}^{k} x^{(j)}, \quad
h=\frac{1}{n+1-k}\sum_{j=k+1}^{n+1} x^{(j)}
\end{equation}
mean that  $(n+1)c=kg+(n+1-k)h,$ �.\,�.
\begin{equation}\label{c_in_gh}
c=\frac{k}{n+1}g+\frac{n+1-k}{n+1}h.
\end{equation}

For proving  (\ref{max_point_y}), it is sufficient to indicate a linear
polynomial 
$p$ which takes in~the~nodes values $\pm 1$ and such that
$p(y)=\|P\|_B$. The equalities $p\left(x^{(j)}\right)=\pm 1$ imply
$$p(y)= \sum_{j=1}^{n+1}p\left(x^{(j)}\right)\lambda_j(y)\leq \lambda(y)=
 \sum_{j=1}^{n+1}|\lambda_j(y)|\leq
\max_{x\in B}\sum_{j=1}^{n+1}
|\lambda_j(x)|=
\|P\|_B.
$$
If $p(y)=\|P\|_B$, then all the values in this chain coincide, therefore,
$\lambda(y)=\|P\|_B.$

Let us show that the above property is fulfilled for the polynomial $p\in 
\Pi_1({\mathbb R}^n)$ with values
\begin{equation}\label{values_of_p}
p\left(x^{(1)}\right)=\ldots=  
p\left(x^{(k)}\right)=-1, \quad 
p\left(x^{(k+1)}\right)=\ldots=  
p\left(x^{(n+1)}\right)=1.
\end{equation}
Since $p$ is a linear function, from
(\ref{c_g_h}) and (\ref{values_of_p}) it follows that
$p(g)=-1,$ $p(h)=1.$ Making use of  (\ref{c_in_gh}), we get
$$p(c)=\frac{k}{n+1}p(g)+\frac{n+1-k}{n+1}p(h)=\frac{n+1-2k}{n+1}.$$
The center of gravity of a regular inscribed simplex coincides with the center of~the~ball. The increment of $p$ is proportional to the distance between the points, hence 
$$\frac{\|g-h\|}{p(g)-p(h)}=\frac{R}{p(y)-p(c)},$$
where $R$ is the radius of the ball.
 Utilizing the found values, we obtain 
\begin{equation}\label{p_y_general}
p(y)=\frac{n+1-2k}{n+1}+\frac{2R}{\|g-h\|}.
\end{equation}
The value of the  latter fraction does not depend on choice of a ball and a regular inscribed simplex. Let us calculate this value for the concrete $S$ and $B$.

Namely, as $S$ we take the regular simplex
with vertices
$$x^{(1)}=e_1, \ \  \ldots, \ \ x^{(n)}=e_n, \quad
x^{(n+1)}=
\left(\frac{1-\sqrt{n+1}}{n},\ldots,\frac{1-\sqrt{n+1}}{n}\right).
$$
The length of  any edge of $S$ is equal to $\sqrt{2}$. This simplex is inscribed
into the ball 
$B=B(x^{h(0)};R)$, where 
$$x^{(0)}=\left(\frac{1-\sqrt{\frac{1}{n+1}}}{n},\ldots,
\frac{1-\sqrt{\frac{1}{n+1}}}{n}\right), 
\quad R=\sqrt{\frac{n}{n+1}}.$$ 
In accordance to  (\ref{c_g_h}), the coordinates of $g$ and $h$ are
$$g_1=\ldots=g_k=\frac{1}{k}, \quad g_{k+1}=\ldots=g_n=0,$$
$$h_1=\ldots=h_k=\frac{1}{n+1-k}\cdot\frac{1-\sqrt{n+1}}{n},$$
$$h_{k+1}=\ldots=h_{n+1}=\frac{1}{n+1-k}\cdot\left(1+\frac{1-\sqrt{n+1}}{n}\right).$$ 
From this,
$$\|g-h\|^2=\left(\frac{1}{k}-\frac{1-\sqrt{n+1}}{n(n+1-k)}\right)^2\cdot k
+\left(\frac{n+1-\sqrt{n+1}}{n(n+1-k)}\right)^2\cdot(n-k).$$
The simple calculation yeilds
$$\|g-h\|^2=\frac{n+1}{k(n+1-k)}.$$ 
Thus, in this case
$$\frac{2R}{\|g-h\|}=2\sqrt{\frac{n}{n+1}}\cdot\frac{\left(k(n+1-k)\right)^{\frac{1}{2}}}{\sqrt{n+1}}=
\frac{2\sqrt{n}}{n+1}\cdot\left(k(n+1-k)\right)^{\frac{1}{2}}.$$

Continuing (\ref{p_y_general}), we can write
$$p(y)=\frac{n+1-2k}{n+1}+\frac{2R}{\|g-h\|}=1-\frac{2k}{n+1}+
\frac{2\sqrt{n}}{n+1}\cdot\left(k(n+1-k)\right)^{\frac{1}{2}}.$$
If $1\leq k\leq \frac{n+1}{2}$, then  the last expression coincides with
$\psi(k)$. The noted inequality is true. Moreover, $k$ coincides with that  number $a$ or $a+1$ on which $\psi(t)$ takes a~bigger value. 
Therefore,
$p(y)= \max\{\psi(a),\psi(a+1)\}=
\|P\|_B.$

The proof is complete. \hfill$\Box$

\smallskip
{\bf Theorem 2.}
{\it In  denotations  of the previous theorem,  $[g,h]$ is the segment of~maximal length in $S$ parallel to vector $gh$. }

\smallskip
{\it Proof.} 
In \cite{nev_fpm_2013}, the author obtained the calculation formulae for  length
and endpoints of the maximal segment in $S$ of a given direction. One can apply these formulas for some simplex 
and take into account the similarity arguments. 
But it is much simpler to use the following characterization of the maximal segment proved in\cite{nev_fpm_2013}
(see there Lemmas ~1 and~2). 
{\it A segment in a simplex parallel to a given vector has maximal length iff every $(n-1)$-dimensional face of the simplex contains at least one endpoint of this segment.}

Let the notation $x=\{\beta_1,\ldots,\beta_{n+1}\}$ means that a point
 $x$ has barycentric coordinates $\beta_1,\ldots,\beta_{n+1}$
with respect to  $S$. 
By $G_j$ denote
$(n-1)$-dimensional face of the simplex not containing the $j$th vertex. For  points of $G_j$, all barycentric coordinates are nonnegative and $\beta_j=0$.
We have
$$
g=\frac{1}{k}\sum_{j=1}^{k} x^{(j)}=\left\{                 
\frac{1}{k},
\ldots, 
\frac{1}{k},
0,\ldots,0\right\},$$
$$
h=\frac{1}{n+1-k}\sum_{j=k+1}^{n+1} x^{(j)}=\left\{
0,\ldots,0,
\frac{1}{n+1-k},
\ldots, 
\frac{1}{n+1-k}\right\}
$$
The number of nonzero barycentric coordinates
in these equalities is equal to 
$k$ and $n+1-k$  respectively. 
Clearly, $g\in G_{k+1},\ldots, G_{n+1}$,
$h\in G_{1},\ldots, G_{k}$. So, every $(n-1)$-dimensional face of  $S$
contains an endpoint of the segment $[g,h]$. Consequently, this segment has maximal length of all the segments of given direction in $S$. Note that this argument is suitable for
any simplex and $k=1,\ldots,n$.
\hfill$\Box$

\section{The Projector's Norm Invariance \\
under an Affine Transform}\label{nev_s3}
In 1948,  F. John \cite{john_1948} proved that any convex body in ${\mathbb R}^n$
contains a unique ellipsoid of maximum volume. Also he gave characterization of those convex bodies for which a maximal ellipsoid is the unit Euclidean ball
$B_n$ (in details see, e.\,g.,\cite{ball_1990}, \cite{ball_1997}).  John's theorem implies the analogous statement which characterizes a unique minimum volume ellipsoid containing a given convex body.

We shall consider a minimum volume ellipsoid containing a given nondegenerate simplex. For brevity, such an ellipsoid will be called {\ it a minimal ellipsoid}. Obviously, a minimal ellipsoid of a simplex is circumscribed around this simplex.
The~center of the ellipsoid coincides with the center of gravity of the simplex.
A~minimal ellipsoid of a simplex is a Euclidean ball iff this simplex is regular. This is equivalent to the well-known fact that the volume of a simplex contained in a ball is maximal iff this simplex is  regular and inscribed into the ball 
(see, e.\,g., \cite{fejes_tot_1964},
\cite{slepian_1969},
\cite {vandev_1992}).

By definition, put $\varkappa_n:=\vo(B_n)$. Denote by
$\sigma_n$ the volume of a regular simplex inscribed into the unit ball $B_n$. 
Suppose $S$ is an arbitrary  $n$-dimensional simplex and  $E$ is the minimal ellipsoid of $S$. 
If a nondegenerate affine transform maps  $S$ into a regular simplex inscribed
into $B_n$, then the image of $E$ under this transform coincides with $B_n$. Hence,
$$
\frac{\vo(E)}{\vo(S)}=\frac{\varkappa_n}{\sigma_n}.$$
It is known that
$$
\varkappa_n=\frac{\pi^
{\frac{n}{2}}}
{\Gamma\left(\frac{n}{2}+1\right)},\qquad
\sigma_n=\frac{1}{n!}\sqrt{n+1}\left(\frac{n+1}{n}\right)^{\frac{n}{2}},
$$
$$
\varkappa_{2m}=\frac{\pi^{m}}{m!},\qquad
\varkappa_{2m+1}=\frac{2^{m+1}\pi^{m}}{(2m+1)!!}=
\frac{2(m!)(4\pi)^m}{(2m+1)!}
$$
(see, e.\,g., \cite{fiht_2001}, \cite{prudnikov_1981},
\cite{nevskii_monograph}). Therefore,   
$$\vo(E)=K_n\, \vo(S), \qquad
K_n:=\frac{\varkappa_n}{\sigma_n}=
\frac{n!\, (\pi n)^{\frac{n}{2}} }{ \Gamma\left(\frac{n}{2}+1\right) (n+1)^{\frac{n+1}{2}} }.
$$
In addition,
$$K_{2m}=\frac{(2m)!(2\pi m)^m}{m!(2m+1)^{m+\frac{1}{2}} },
\qquad
K_{2m+1}=2^{m+\frac{1}{2}}\left(2-\frac{1}{m+1}\right)^{m+\frac{1}{2}} \pi^m m!.$$

The value $K_n$  is included in the lower bound of the norm of a projector with nodes in $B_n$.  
Let $\chi_n(t)$ be {\it  the standardized Legendre polynomial of degree $n$}: 
$$\chi_n(t):=\frac{1}{2^nn!}\left[ (t^2-1)^n \right] ^{(n)}.$$
There exists a constant $C>0$ not depending on $n$ such that for any
interpolation projector $P:C(B_n)\to\Pi_1({\mathbb R}^n)$
\begin{equation}\label{est_norm_P_via_chi_n}
\|P\|_{B_n}\geq \chi_n^{-1}(K_n)>C\sqrt{n}.
\end{equation}
Inequalities \eqref{est_norm_P_via_chi_n} were obtained by the author in
 \cite{nevskii_mais_2019_3}. The right-hand estimate holds true,   if we take, e.\,g., 
$$C=
\frac{  \sqrt[3]{\pi} }{\sqrt{12e}\cdot\sqrt[6]{3} } =  0.2135...$$

Assume $S$ and $S^\prime $ are nondegenerate simplices in ${\mathbb R}^n$ with vertices   
$x^{(j)}, \ldots, x^{(n+1)}$ and 
$y^{(1)}, \ldots, y^{(n+1)}$ respectively. Let $\bf S$ be the vertex matrix of  $S$. Denote  by
${\bf Y}$  the~$n\times(n+1)$-matrix whose $j$th column contains the coordinates of $y^{(j)}$. Let $\lambda_1,\ldots,\lambda_{n+1}$ be the basic Lagrange polynomials of $S$.

\smallskip
{\bf Lemma 1.} {\it There exists a unique affine transform
$F$ of  space ${\mathbb R}^n$ which maps
$S$ into $S^\prime$ and such that $y^{(j)}=F\left(x^{(j)}\right)$.  The equality $y=F(x)$
is equivalent to any relation
\begin{equation}\label{theor_affine_transform_matrices}
\left(
\begin{array}{c}
y_1\\
\vdots\\
y_{n}
\end{array} 
\right)=
{\bf Y}\left({\bf S}^{-1}\right)^{T}
\left(
\begin{array}{c}
x_1\\
\vdots\\
x_{n}\\
1
\end{array} 
\right),
\end{equation}
\begin{equation}\label{theor_affine_transform_Lagrange_polynomials}
y=\sum\limits_{j=1}^{n+1} \lambda_j(x) y^{(j)}.
\end{equation}
}

\smallskip
{\it Proof.}
Each nondegenerate affine transform of  ${\mathbb R}^n$ has the form $F(x)=A(x)+b$, where $A:{\mathbb R}^n\to {\mathbb R}^n$  is a nondegenerate
linear operator. Let  ${\bf A}=(a_{ij})$ be the~matrix of the operator $A$ in the canonical basis. In coordinate form, the equality $y=A(x)+b$  is equivalent to the relation 
$$
\left(
\begin{array}{c}
y_1\\
\vdots\\
y_{n}
\end{array} 
\right)=
\left(
\begin{array}{cccc}
a_{11}&\ldots&a_{1n}&b_1\\
\vdots&\vdots&\vdots&\vdots\\
a_{n1}&\ldots&a_{nn}&b_n
\end{array}
\right) 
\left(
\begin{array}{c}
x_1\\
\vdots\\
x_{n}\\
1
\end{array} 
\right).
$$
Define ${\bf M}$ as the 
$n\times(n+1)$-matrix standing on right-hand side. 
The conditions $y^{(j)}=F\left(x^{(j)}\right)$  are equivalent to the equality ${\bf Y}={\bf M}{\bf S}^{T}$.
Cosequently,
 $${\bf M}={\bf Y}\left({\bf S}^{T}\right)^{-1}={\bf Y}\left({\bf S}^{-1}\right)^{T}.$$
This means that an affine transform  satisfying the conditions of the theorem is~unique and has the form 
\eqref{theor_affine_transform_matrices}. 

Since $\lambda_j\in \Pi_1({\mathbb R}^n)$ and $\lambda_j\left(x^{(k)}\right)$ $=$ ~$\delta_j^k$, the equality \eqref{theor_affine_transform_Lagrange_polynomials} also gives an affine transform 
$y=F(x)$ such that $F\left(x^{(k)}\right)=y^{(k)}$. From the uniqueness of $F$, 
it follows that \eqref{theor_affine_transform_Lagrange_polynomials}  is equivalent to  \eqref{theor_affine_transform_matrices}. This equivalence  can be proved also directly.  Let us rewrite \eqref{theor_affine_transform_Lagrange_polynomials} in the coordinate form using the coefficients of the polynomials $\lambda_j$:
$$y=\sum\limits_{j=1}^{n+1} \lambda_j(x) y^{(j)}=
\sum\limits_{j=1}^{n+1} \left(\sum_{k=1}^n l_{kj}x_k+l_{n+1,j}\right)y^{(j)},$$ 
$$y_i=\sum\limits_{j=1}^{n+1} \left(\sum_{k=1}^n l_{kj}x_k+l_{n+1,j}\right)y_i^{(j)}=
\sum\limits_{j=1}^{n+1} \left(\sum_{k=1}^n l_{kj}x_k\right)y_i^{(j)}+
\sum\limits_{j=1}^{n+1} l_{n+1,j}y_i^{(j)}=$$
$$=\sum_{k=1}^n\left(\sum_{j=1}^{n+1}  y_i^{(j)} l_{kj}\right)x_k+
\sum_{j=1}^{n+1}y_i^{(j)}l_{n+1,j}.$$
Thus, \eqref{theor_affine_transform_Lagrange_polynomials} means that 
$$y_i=
\sum_{k=1}^n a_{ik}x_k+b_i,\qquad 
a_{ik}=\sum_{j=1}^{n+1}  y_i^{(j)} l_{kj}, \quad  b_i=\sum_{j=1}^{n+1}y_i^{(j)}l_{n+1,j}.
$$
Since  ${\bf S}^{-1}$ $=(l_{ij})$, these equalities are equivalent to \eqref{theor_affine_transform_matrices}.
\hfill$\Box$


\smallskip
We note that the norm of an interpolation projector is invariant under an affine transform.

\smallskip
{\bf Theorem 3.} {\it
  Suppose $\Omega$ is a convex body in ${\mathbb R}^n$ containing
 a nondegenerate simplex $S$, $\Omega^\prime$ and
$S^\prime$ are their images under a nondegenerate affine transform, \linebreak
$P:C(\Omega)\to\Pi_1({\mathbb R}^n)$ and
$P^\prime:C(\Omega^\prime)\to\Pi_1({\mathbb R}^n)$ are interpolation projectors
with the nodes in the vertices of $S$  and $S^\prime$ respectively.
Then $\|P\|_\Omega=\|P^\prime\|_{\Omega^\prime}.$
}

\smallskip
{\it Proof.} 
Let $x_1,\ldots,x_{n+1}$ be the vertices of the simplex $S$. We will assume that
the verices of~the simplex $S^\prime$ are numerated so that  
$y^{(j)}=F\left(x^{(j)}\right)$. Under this condition, the set of barycentric coordinates of an arbitrary point $x\in {\mathbb R}^n$ with respect to $S$ coincides with the
set of barycentric coordinates of the point 
$y=F(x)$ with respect to  $S^\prime$. This follows from the equalities
$$x=\sum\limits_{j=1}^{n+1} \lambda_j(x) x^{(j)}, \quad y=\sum\limits_{j=1}^{n+1} \lambda_j(x) y^{(j)}.$$
The second equality coincides with  relation
 \eqref{theor_affine_transform_Lagrange_polynomials} of Lemma 1.
In accordance with formula  \eqref{norm_P_Omega_by_barycentric_coord},
 we have 
$$
\|P\|_\Omega=\max \left \{ \sum_{j=1}^{n+1} |\beta_j|:  \ \sum_{j=1}^{n+1} \beta_jx^{(j)} \in \Omega, \ \sum_{j=1}^{n+1} \beta_j=1\right\}=$$
$$=\max \left \{ \sum_{j=1}^{n+1} |\beta_j|:  \ \sum_{j=1}^{n+1} \beta_jy^{(j)} \in \Omega^\prime, \ \sum_{j=1}^{n+1} \beta_j=1\right\}=\|P^\prime\|_{\Omega^\prime}.
$$
\hfill$\Box$

\smallskip
{\bf Corollary 1.} 
{\it Suppose $S$ is a nondegenerate  
simplex with minimal ellipsoid $E$ and $S^\prime$ is an arbitrary regular simplex inscribed into $B_n$. 
If $P:C(E)\to\Pi_1({\mathbb R}^n)$ and 
$P^\prime:C(B_n)\to\Pi_1({\mathbb R}^n)$ are the projectors having nodes in the vertices of $S$ and $S^\prime$ respectively, then
$\|P\|_E=\|P^\prime\|_{B_n}.$
}


\smallskip
{\it Proof.} Consider the nondegenerate affine transform which maps the simplex
$S$ into the regular simplex  $S^\prime$. This transform maps the ellipsoid $E$ into the ball $B_n$. It remains to apply Theorem 3 in the case when $\Omega$ is the minimal ellipsoid of $S$.
  \hfill$\Box$

\smallskip
Let us supplement Corollary 1 with the following remark. Denote here by $\lambda_j$
the basic Lagrange polynomials  of a simplex $S$.
The maximum points of the function   $\lambda (x) $ = $ \sum | \lambda_j (x) | $ lying
in the minimal ellipsoid $E$  have the same geometric description that is formulated in Theorem 1. In the condition of this theorem, we must replace the regular simplex   by an arbitrary one, and the circumscribed ball by the minimal ellipsoid of the simplex. At the specified points of the border \linebreak of the ellipsoid, $\lambda (x) $ takes maximal value  equal to $ \| P \|_E $. This result can be established according to the scheme above. 

\smallskip
{\bf Corollary 2.} {\it There exists a universal constant $C>0$ such that for every ellipsoid $E\subset{\mathbb R}^n$ and every interpolation projector having the
nodes in $E$ we have
$\|P\|_{E}\geq \chi_n^{-1}(K_n)>C \sqrt{n}.$}

\smallskip
This follows immediately from \eqref{est_norm_P_via_chi_n} and Corollary 1.

\section{On Some Extremal Property of a Regular Simplex Inscribed into a Ball}\label{nev_s4}

Consider  a nondegenerate simplex $S\subset {\mathbb R}^n$. Let $E$ be the minimal ellipsoid containing  $S$. Fix a natural number $m\leq\frac{n}{2}$. To each set of $ m $ vertices of $S$ assign the point $ y \in E $ defined as follows. Let $ g $ be
the center of gravity of the $ (m-1) $-dimensional face of $ S $ containing the selected vertices, and let $ h $ be the center of~gravity of the $(n-m)$-dimensional face containing the remaining $ n + 1- m $ vertices. Then $ y $ is the intersection point of the straight line $ (gh) $ with the boundary of $E$ in~the~direction from $ g $ to $ h $.

  Now we formulate the following conjecture.
  
\smallskip
(H1) {\it For a given  $m\leq\frac{n}{2}$ and any nondegenerate simplex  $S\subset B_n$,
there exists a set of  $m$  vertices of $S$ such that $y\in B_n$.}

\smallskip
A stronger version of the hypothesis asserts that the specified property holds \linebreak
{\it for any $ m  \leq \frac {n} {2} $} (H2).
For our purposes, it is sufficient that (H1) was true
for $ m = k(n) $. The number $ k = k(n) $ is defined in Section 2.

\smallskip
{\bf Theorem 4.} {\it For $m=1$ conjecture (H1) holds true.}

\smallskip
{\it Proof.}
Suppose $S$ is a simplex with vertices $x^{(j)}\in B_n$ and the center of
gravity  
$c$. The center of the minimal ellipsoid containing $S$ also lies in $c$. 
Hence, in the~case $m=1$ the points $y$ has the form
$y^{(j)}=2c- x^{(j)},$ $j=1,\ldots,n+1.$
 We need to~show that there exists a vertex $x$ of the simplex such that $\|2c-x\|\leq 1$. 
 Since $S$ is~nondegenerate, for some vertex $x$ we have 
 $(c,x-c)\geq 0$. This means that
$$\|2c-x\|^2=(2c-x,2c-x)=
4(c,c-x)+\|x\|^2\leq \|x\|^2\leq 1,$$
i.\,e., the vertex $x$ is suitable.
The theorem is proved.
\hfill$\Box$

\smallskip
Denote by $\theta_n(B_n)$ the minimal norm of an interpolation projector 
$P:C(B_n)\to \Pi_1({\mathbb R}^n)$ with the nodes in $B_n$. By $P^\prime$ denote a projector
whose nodes coincide with the vertices of a regular simplex $S^\prime$ inscribed into $B_n$.

\smallskip
{\bf Theorem 5.}
 {\it 
 Suppose (H1) is true for   $m=k(n)$.  Then
$\theta_n(B_n)=\|P^\prime\|_{B_n}.$
}

\smallskip
{\it Proof.}
Consider an arbitrary projector
$P$ with the nodes $x^{(j)}\in B_n$. Let $S$ be the simplex with these vertices and let  $\lambda_j$ be the basic Lagrange polynomials corresponding to $S$.
Denote by $E$ the minimal ellipsoid of the simplex. 
Since $S\subset B_n$, for some set consisting of $k=k(n)$ vertices of the simplex the corresponding
point $y$ is contained in the ball.  
Let us fix $y$ and write the following relations:
$$\|P^\prime\|_{B_n}=\|P\|_E=\max_{x\in E}\sum_{j=1}^{n+1} |\lambda_j(x)|=
\sum_{j=1}^{n+1} |\lambda_j(y)|\leq \max_{x\in B_n}\sum_{j=1}^{n+1} |\lambda_j(x)|=
\|P\|_{B_n}.$$
We made use of the formula for the projector norm, Theorem 1, Corollary 1, and~also the remark after this 
corollary. The inequality in the above chain follows from the~condition  $y\in B_n$. Note that if
 $y$ lies inside the ball, then this equality becomes strict. 
 
 Therefore, for any projector with nodes in  $B_n$ we have
 $\|P^\prime\|_{B_n}\leq \|P\|_{B_n}$. This implies that $\theta_n(B_n)=\|P^\prime\|_{B_n}$.
 The proof is complete. 
\hfill$\Box$

\smallskip
{\bf Corollary 3.} {\it
If $1\leq n \leq 4$, then $\theta_n(B_n)=\|P^\prime\|_{B_n}.$
}

\smallskip
{\it Proof.}
In the  case $n=1$, the proposition is equivalent to the fact that the~norm of~an~interpolation projector $P:C[-1,1]\to \Pi_1({\mathbb R})$ becomes minimal  for the projector having the nodes at the endpoints of the segment $[-1,1]$. If $2\leq n \leq 4$, then $k(n)=1$, and the required result follows immediately from Theorems 4 and 5.  
\ \hfill$\Box$

\smallskip
Corollary 3 was proved in  \cite{nev_ukh_mais_2019_2} by another method
suitable only for dimensions $n$ with the property $k(n)=1$. 
However, starting from $n=5$, we have \linebreak $k(n)>1$ (see~\cite{nev_ukh_mais_2019_2}).
Nethertheless,  the equality $\theta_n(B_n)=\|P^\prime\|_{B_n}$ still can be obtained on~the~way directed by Theorem 5.

In the propositions of this section, the unit ball $B_n$ may be replaced by~an~arbitrary Euclidean ball $B$;
this leads to the equivalent results.


\end{document}